\documentclass{elsart}

\usepackage{graphicx}
\usepackage{mathrsfs}
\usepackage{amssymb}
\usepackage{amsmath}
\usepackage{graphics}
\usepackage{latexsym}
\usepackage{amsfonts}
\usepackage{color}
\usepackage{ctable}
\usepackage{verbatim}

\numberwithin{equation}{section}

\newcommand{\be}{\begin{eqnarray*}}
\newcommand{\ee}{\end{eqnarray*}}
\newcommand{\ben}{\begin{equation}}
\newcommand{\een}{\end{equation}}
\newcommand{\lb}[1]{\left[\begin{array}{#1}}
\newcommand{\rb}{\end{array}\right]}
\newcommand{\lp}[1]{\left(\begin{array}{#1}}
\newcommand{\rp}{\end{array}\right)}

\newcommand{\leftd}[1]{\left\{\begin{array}{#1}}
\newcommand{\rightd}{\end{array}\right.}

\def\A {\mathbf{A}}
\def\B {\mathbf{B}}
\def\C {\mathbf{C}}

\def\H {\mathbf{H}}
\def\I {\mathbf{I}}

\def\L {\mathbf{L}}
\def\M {\mathbf{M}}

\def\P {\mathbf{P}}

\def\R {\mathbf{R}}

\def\T {\mathbf{T}}

\def\b {\boldsymbol{b}}

\def\e {\boldsymbol{e}}

\def\y {\boldsymbol{y}}

\def\bbeta {\boldsymbol{\beta}}

\def\Nb {\mathbb{N}}
\def\Rb {\mathbb{R}}

\begin{document}
\begin{frontmatter}

\title{Strong noise estimation in cubic splines}
\author[South]{A. Dermoune}
\ead{Azzouz.Dermoune@math.univ-lille1.fr},
\author[Oxford]{A. El Kaabouchi}
\ead{aek@ismans.fr}

\address[South]{Laboratoire Paul Painlev\'e, USTL-UMR-CNRS 8524. UFR de Math\'ematiques, B\^at. M2, 59655 Villeneuve d'Ascq Cedex, France} 
\address[Oxford]{Institut Sup\'erieur des Mat\'eriaux et M\'ecaniques Avanc\'es, 44, Avenue Bartholdi, 72000, Le Mans, France}

\begin{abstract} 

The data $(y_i,x_i)\in \Rb\times [a,b], i=1,\ldots,n$ satisfy
$y_i=s(x_i)+e_i$ where $s$ belongs to  the set of cubic splines. The unknown noises $(e_i)$ are such that  $var(e_I)=1$
for some $I\in \{1, \ldots, n\}$ and $var(e_i)=\sigma^2$ for $i\neq I$. 
We suppose that the most important noise is $e_I$, i.e. the ratio $r_I=\frac{1}{\sigma^2}$ is larger than one.  
If the ratio $r_I$ is large, then we show, for all smoothing parameter, that the penalized least squares estimator of the $B$-spline basis 
recovers exactly the position $I$ and the sign of the most important noise $e_I$.   

\end{abstract}

\end{frontmatter}

\section{Linear inverse problem: General setting} 
The data $(y_i,x_i)\in \Rb\times [a,b], i=1, . . ., n$ satisfy
\begin{equation}
y_i=s(x_i)+e_i.
\end{equation}
The map $s: [a, b]\to \Rb$ is unknown, and $(e_i)$ are the error of measurements, also called the noise and is unknown. 
We suppose that $s$ belongs to a set ${\cal C}$ of functions, and we are interested in the estimation of $s\in {\cal C}$ using the data $(y_i,x_i), i=1, \ldots, n$.
Suppose that ${\cal C}$ has a basis $(b_j)_{j=1, \ldots, d}$. In this case, each map $s\in {\cal C}$
is determined by its coordinates $\bbeta=(\beta_1 \ldots \beta_d)^T$ in the latter basis, i.e., 
$\forall x\in [a,b],\quad s(x)=\sum\limits_{j=1}^d\beta_j b_j(x)$. 
Hence, for each $i$, $s(x_i)=\b(x_i)\bbeta$ with $\b(x_i)=\left(b_1(x_i) \ldots b_d(x_i)\right)$ is the $(1,d)$ matrix. 
If we introduce the $(n,d)$ matrix 
\ben
\B=\lp{c}\b(x_1)\\ \vdots\\ \b(x_n)\rp,
\een
then, the data $\y=(y_1 \ldots y_n)^T$ and the noise $\e=(e_1 \ldots e_n)^T$ satisfy the linear system 
\ben 
\y=\B\bbeta+\e.
\label{LIP}
\een 
The latter is known as the linear regression in Statistic community and the linear inverse problem 
in the Inverse problem community. This problem is ill-posed, because the transformation $\bbeta\mapsto \B\bbeta$ is not invertible. 
Moreover, the noise $\e$ is not known. 
  
One way to estimate the parameter $\bbeta$ and the noise $\e$ is to use 
the generalized penalized least square estimators. It works as following. 
We fix a matrix $\M$ having $n$ columns, and we consider 
the ellipsoide quasi-norm $\|\cdot\|_\M$ defined by $\|x\|_\M^2=x^T\M^T\M x$. 
We propose, for each $\lambda >0$ and for each matrix $\L$ having $d$ columns, the minimizers 
\ben
\hat{\bbeta}(\lambda,\M,\L)\in \arg\min_{\bbeta}\{\|\y-\B\bbeta\|_\M^2+\lambda \|\L\bbeta\|^2\}
\een
as an estimator of the vector $\bbeta$. 
The quantity $\|\y-\B\bbeta\|_\M^2$ is the square of the residual error with respect to 
the metric defined by the quasi-norm $\|\cdot\|_\M$, and $\|\L\bbeta\|^2$ is called the penality. 
The parameter $\lambda$ is called the smoothing parameter. We have easily the following results.   
\begin{prop} 
The set of the minimizers of the latter optimization is given 
by the following normal equation 
\ben 
(\B^T\M^T\M\B+\lambda \L^T\L)\bbeta=\B^T\M^T\M\y.
\een
If $N(\M\B)\cap N(\L)=\{0\}$, then $\hat{\bbeta}(\lambda,\M,\L)$ is unique and is given by 
\ben 
\hat{\bbeta}(\lambda,\M,\L)=(\B^T\M^T\M\B+\lambda \L^T\L)^{-1}\B^T\M^T\M\y
:=\H(\lambda,\M,\L)\y.
\label{hatmatrix}
\een
Here, $N(\A)$ denotes the null-space of the matrix $\A$. 
\end{prop} 
The minimizer $\H(\lambda,\M,\L)\y$
is proposed as an estimator of the parameter $\bbeta$. Hence, $\B\H(\lambda,\M,\L)\y$ is an estimator 
of $\B\bbeta$ and $\y-\B\H(\lambda,\M,\L)\y$ is an estimator of the noise $\e$.  
The map $x\in [a,b]\mapsto \sum\limits_{j=1}^d\H(\lambda,\M,\L)\y(j)b_j(x)$ is an estimator of the map $s$. 
The performance of these estimators depends clearly on the matrix $\H(\lambda, \M,\L)$ known as the hat matrix in Statistics community.
\begin{prop}\label{elden}  The limit   
$\lim_{\lambda\mapsto 0^+}\H(\lambda,\M,\L):=H(0,\M,\L)$
exists and is equal to the $\M\L$-weighted pseudoinverse of $\B$ defined by 
\ben 
\B_{\M\L}^+=(\I-(\L\P_{\M\B})^+\L)(\M\B)^+\M.
\een 
Here $\A^+$ denotes the Moore-Penrose inverse of $\A$ and 
\ben 
\P_\A=\I-\A^+\A
\een 
denotes the orthogonal projection on $N(\A)$.

It follows that the estimator   
$\hat{\bbeta}(\lambda,\M,\L)$ converges to $\B_{\M\L}^+\y:=\hat{\bbeta}(0,\M,\L)$ as $\lambda\mapsto 0$. 
Moreover, we can show that 
\ben 
\hat{\bbeta}(0,\M,\L)=\arg\min_{\bbeta}\{\|\L\bbeta\|^2:\quad\bbeta\in \arg\min\|\y-\B\bbeta\|_\M^2\}.
\een 
In particular, if $\B$ has maximal rank, i.e., $R(\B)=\R^n$, then  
\ben 
\hat{\bbeta}(0,\M,\L)=\arg\min_{\bbeta}\{\|\L\bbeta\|^2:\quad \y=\B\bbeta\},
\een 
or equivalently 
\ben
\B\B_{\M\L}^+\B=\B. 
\label{exactnoisezero}
\een
Observe that, if $\M$ is invertible, then  
\ben 
\B_{\M\L}^+\B=\B_{\L}^+\B. 
\een
On the other hand, the limit 
\ben 
\lim_{\lambda\mapsto +\infty}(\B^\T\M^T\M\B+\lambda \L^T\L)^{-1}\B^T\M^T\M
\een 
exists and is equal to 
\ben 
\C_{\M\B\L}:=(\M\B\P_\L)^+\M.
\een 
In the classical case $\M=\I$ and $\L=\I$, the matrix  
\ben 
\H(\lambda,\I,\I)=(\B^T\B+\lambda\I)^{-1}\B^T
\een 
is known as the Tikhonov regularization of the Moore-Penrose inverse of the matrix $\B$. 
\end{prop} 

\begin{pf}
It is a consequence of known results\cite{Elden}.
\end{pf}

In the sequel, we  suppose that the noise $\e$ is Gaussian with the covariance matrix $\C=diag(\sigma_i^2)$. 
In this case, the natural choice of $\M$ is the weight matrix $\C^{-1/2}$.
We suppose that the variance $\sigma^2_I=1$ for some $I$ and $\sigma^2_i=\sigma^2$ for all $i\neq I$. 
The set of functions ${\cal C}$ is the set of cubic splines. The true signal is an element 
of the $B$-spline basis. 
We consider, for each $\lambda >0$, the noise estimator $\hat{\e}(\lambda)=\y-\B\H(\lambda,\M,\L)\y$.  
The aim of our work is to show that if $\sigma^2$ is small, then $\hat{\e}(\lambda)$
recovers exactly the position and the sign of the most important noise. 
Section 2 recalls some cubic splines results. In Section 3 we present our numerical results.

\section{Cubic splines} 
Schoenberg introduced in \cite{Schoenberg} the terminology spline for a certain type of piecewise polynomial interpolant. 
The ideas have their roots in the aircraft and shipbuilding industries. 
Since that time, splines have been shown to be applicable and effective for a large number tasks in interpolation and approximation.
Various aspect of splines and their applications can be found in \cite{Micula1999}.  
Let $a=\kappa_0 < \kappa_1 < \ldots< \kappa_{K+1}=b$ be a sequence of increasing real numbers. Spline interpolation can be described as following. 
A map $s$ belongs to the set $S_3(\kappa_0,\ldots, \kappa_{K+1})$ of cubic splines with knots $(\kappa_0 < \kappa_1, \ldots, \kappa_{K+1})$ if 
\ben 
s(x)=p_i+q_i(x-\kappa_i)+\frac{u_i}{2}(x-\kappa_i)^2+\frac{v_i}{6}(x-\kappa_i)^3
\een 
for every $x\in [\kappa_i,\kappa_{i+1})$. 
Let $(b_j:=S_{j,4}:\quad j=-3, \ldots, K)$ be the $B$-spline basis functions of the set $S_3(\kappa_0,\ldots, \kappa_{K+1})$.     

Before going further, let us recall the famous result of \cite{Schoenberg} and \cite{Reinsch}. 
If the data $(y_i,\kappa_i):\,i=0, \ldots, K+1$, then the minimizer of 
\ben
\min_{s\in S_3( \kappa_0,\ldots, \kappa_{K+1})}\{\lambda \int_{a}^{b}|s^{(2)}(x)|^2dx+\sum_{i=0}^{K+1}|s(\kappa_i)-y_i|^2\}
\een 
is the natural cubic spline, i.e. such that $s^{(2)}(\kappa_0)=s^{(2)}(\kappa_{K+1})=0$, where $s^{(2)}$ is the second derivative of $s$.   
Observe that the penalized matrix $\L$ is defined by 
\ben 
\|\L\bbeta\|^2=\int_{a}^{b}|s^{(2)}(x)|^2dx.
\een 

Let us calculate the matrix $\L$. The unknown vector $\bbeta$ belongs to $\Rb^{K+4}$. From the derivative formula for $B$-spline functions \cite{de Boor}, ch. X. we have 
\ben 
\sum_{j=-3}^K\bbeta_jS_{j,4}^{(2)}(x)=\sum_{j=-1}^{K}\bbeta_j^{(2)}S_{j,2}(x),
\een 
where the vector $\bbeta^{(2)}=\Delta_2\bbeta$ and $\Delta_2$ denotes the matrix corresponding to the weighted difference operator. 
If we denote by $\R$ the matrix with entries 
\ben 
R_{ij}=\int_a^b S_{j,2}(x)S_{i,2}(x)dx,\quad i,j=-1,\ldots, K,
\een
then 
\ben 
\int_a^b|S^{(2)}(x)\bbeta|^2dx=\bbeta^t\Delta_2^t\R\Delta_2\bbeta=\|\L\bbeta\|^2, 
\een 
with $\L=\R^{1/2}\Delta_2$. 

In the sequel we suppose that the data $(y_i,x_i),\quad i=1, \ldots, n$ with $(x_i: \quad i=1, \ldots, n)$
do not necessarily coincide with the knots $(\kappa_i:\quad i=0, \ldots, K+1)$.  

We want to study the estimator $\y-\B\H(\lambda,\M,\L)\y$ with respect to the smoothing parameter $\lambda >0$.
More precisely we want to recover the position and the sign of the most important noise. 

\section{Numerical computation}
We consider $a=0$, $b=1$, $K,n\in\Nb^*$ and $(\kappa_i)_{i=0,\ldots,K+1}$ with ${\kappa_{i+1}-\kappa_i=\frac{b-a}{K+1}}$ for all $i\in\{0,\ldots,K\}$. 
The data $(y_i,x_i)_{i=1, \ldots, n}$ are such that
${x_{i+1}-x_i=\frac{b-a}{n-1}}$ for all $i\in\{1,\ldots,n-1\}$. The 
model is $\y=\B\delta_j+\e$. 

The following show that for all $j$ and for each smoothing parameter $\lambda$ the noise estimator $[\I-\B\H(\lambda, \M)]\y$ 
recovers exactly the position $I$ and the sign of the most important noise. We fix the variance $\sigma^2\in (0,1)$, and we consider, for each realization of the noise $\e$, the maps $I(\e)$ 
and $sgn(\e)$ defined respectively by:    
\ben 
\lambda\in (0, 10)\to \arg\max_{i}|[\I-\B\H(\lambda, \M)]\y(i)|=I(\e,\lambda)\in \{1, \ldots, n\},
\een
\ben
\lambda\in (0, 10)\to sgn(\e,\lambda)=sign([\I-\B\H(\lambda, \M)]\y(I)). 
\een 
where $sign(x)=-1$ if $x < 0$, $sign(x)=1$ if $x >0$.

We repeat 100 realizations $(\e^{(k)}:\quad k=1, \ldots, 100)$. We calculate 
the proportion
\ben
p_1(\sigma,\lambda,n)=\frac{1}{100}\sum_{k=1}^{100}1_{[I(\e^{(k)},\lambda)\neq I]},
\een
i.e. the probability that the estimator $[\I-\B\H(\lambda, \M)]\y$ does not recover 
the position $I$ of the strong noise $e_I$. We also calculate the probability that 
the estimator $[\I-\B\H(\lambda, \M)]\y$ does not recover 
the $sign(e_I)$ of the strong noise $e_I$, i.e. 
\ben
p_2(\sigma,\lambda,n)=\frac{1}{100}\sum_{k=1}^{100}1_{[sgn(\e^{(k)},\lambda)\neq sign(e(I))]}. 
\een
The probability that the path $\lambda\in (0,10)\to I(\e^{(k)},\lambda)$ does not coincide with 
the position $I$ of the most important noise $e_I$ is equal to  
\ben
p_3(\sigma,n)=\frac{1}{100}\sum_{k=1}^{100}1_{[I(\e^{(k)})\neq I]}. 
\een
The probability that the path $\lambda\in (0,10)\to sign(\e^{(k)},\lambda)$ does not coincide with 
the sign of the most important noise $e_I$ is equal to  
\ben
p_4(\sigma,n)=\frac{1}{100}\sum_{k=1}^{100}1_{[sgn(\e^{(k)})\neq sign(e_I)]}. 
\een

Below we plot $\lambda\in (0,10)\to p_l(\sigma,\lambda,n)$ for $l=1, 2$ and for fixed $\sigma$, $\sigma\in (0,1.5)\to p_l(\sigma,\lambda,n)$ for $l=1, 2$ and for fixed $\lambda$, and $\sigma\to p_l(\sigma,n)$ for $l=3, 4$.  

\noindent The following example illustrates our results when $K=4$, $j=3$, $I=1$, $n=6$, $n=8$, $n=10$ and $n=20$.

\noindent  $\bullet$ Plots of $p_1(\sigma,\lambda,n)$ and $p_2(\sigma,\lambda,n)$   for fixed $\sigma.$

\begin{center}
\includegraphics{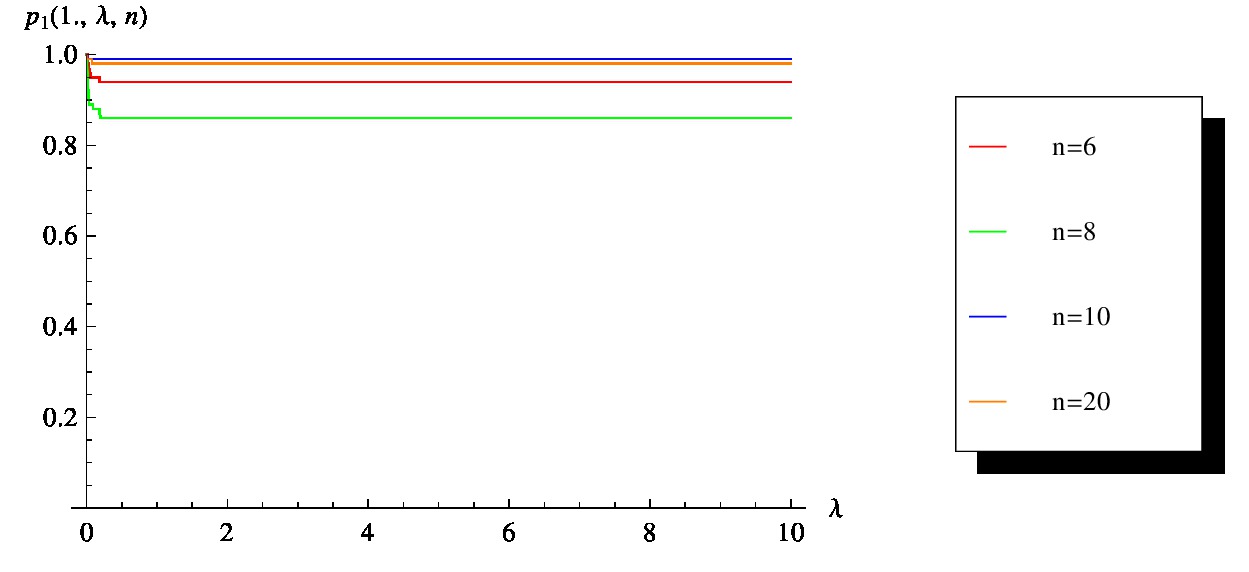}
\end{center}

If the noise is white, then there is no dominating noise 
among $e_1, \ldots, e_n$.  Numerical result are coherent. It tells us that 
the probability that the estimator 
$[\I-\B\H(\lambda, \M)]\y$ recovers the position of the most important noise 
is very small. 

\begin{center}
\includegraphics{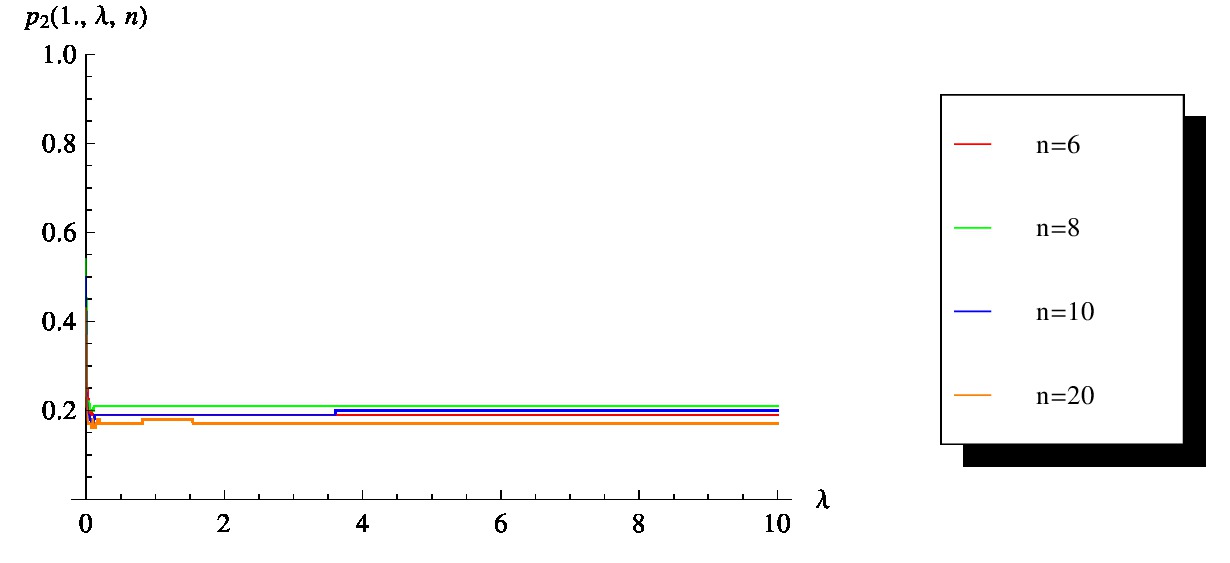}
\end{center}

Numerical result shows that even the noise is white  
the probability that the estimator 
$[\I-\B\H(\lambda, \M)]\y$ recovers the sign of the noise is nearly equal 
to 0.8.  

\begin{center}
\includegraphics{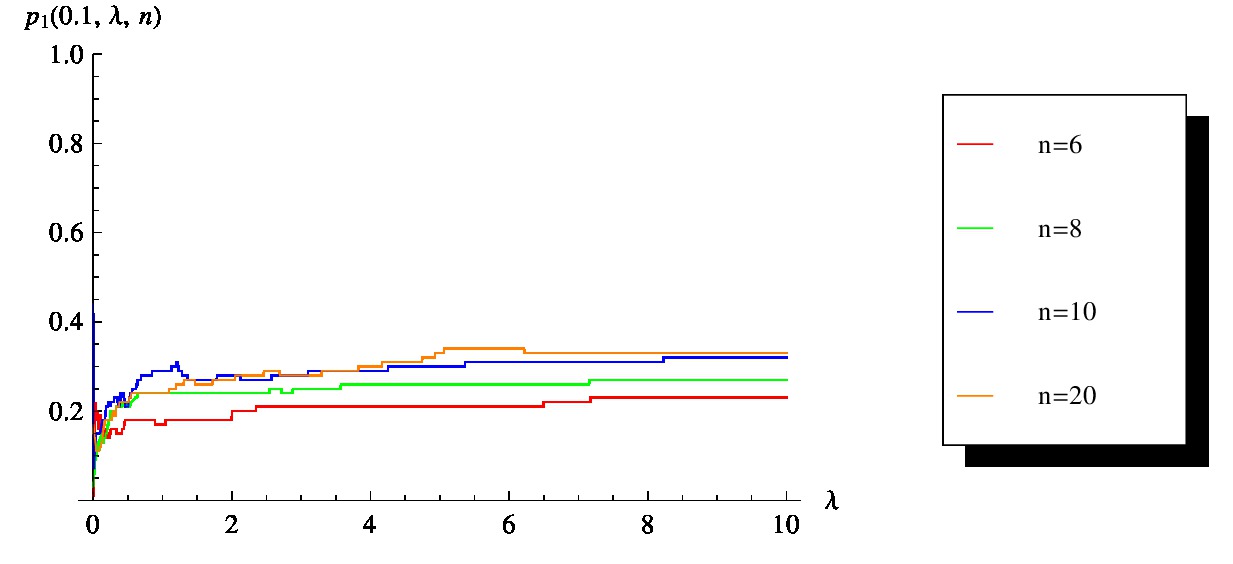}
\end{center}

If the noise has a dominating component, then 
the probability that the estimator 
$[\I-\B\H(\lambda, \M)]\y$ recovers the position of the most important noise 
belongs to $(0.7, 0.9)$ for all $\lambda\in (0,10)$.  

\begin{center}
\includegraphics{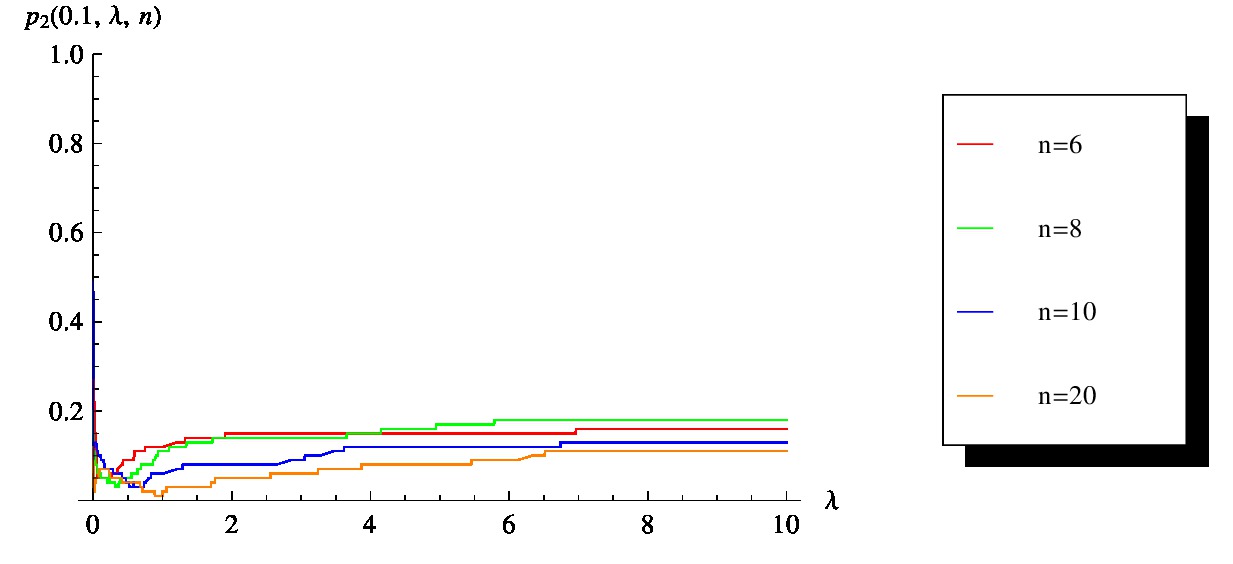}
\end{center}

If the noise has a dominating component, then 
the probability that the estimator 
$[\I-\B\H(\lambda, \M)]\y$ recovers the sign of the most important noise 
belongs to $(0.8, 1)$ for all $\lambda\in (0,10)$.

\noindent  $\bullet$ Plots of $\sigma\to p_1(\sigma,\lambda,n)$ and $\sigma\to p_2(\sigma,\lambda,n)$ for fixed $\lambda$. Numerical results show that both are increasing, but 
$\sigma\to p_1(\sigma,\lambda,n)$ increases quickly than $\sigma\to p_2(\sigma,\lambda,n)$. 

\begin{center}
\includegraphics{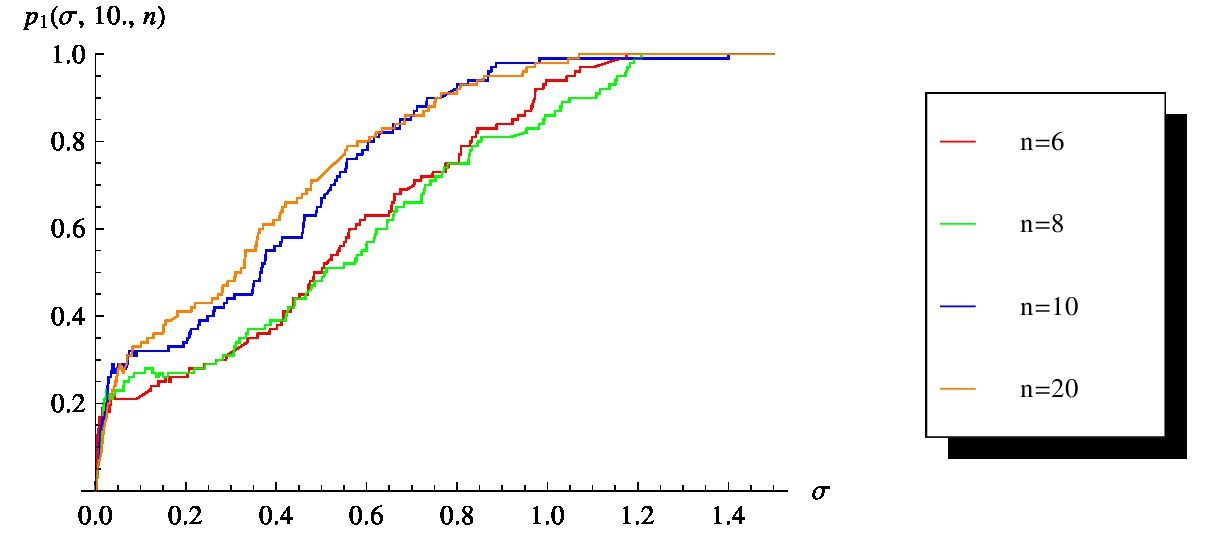}
\end{center}

\begin{center}
\includegraphics{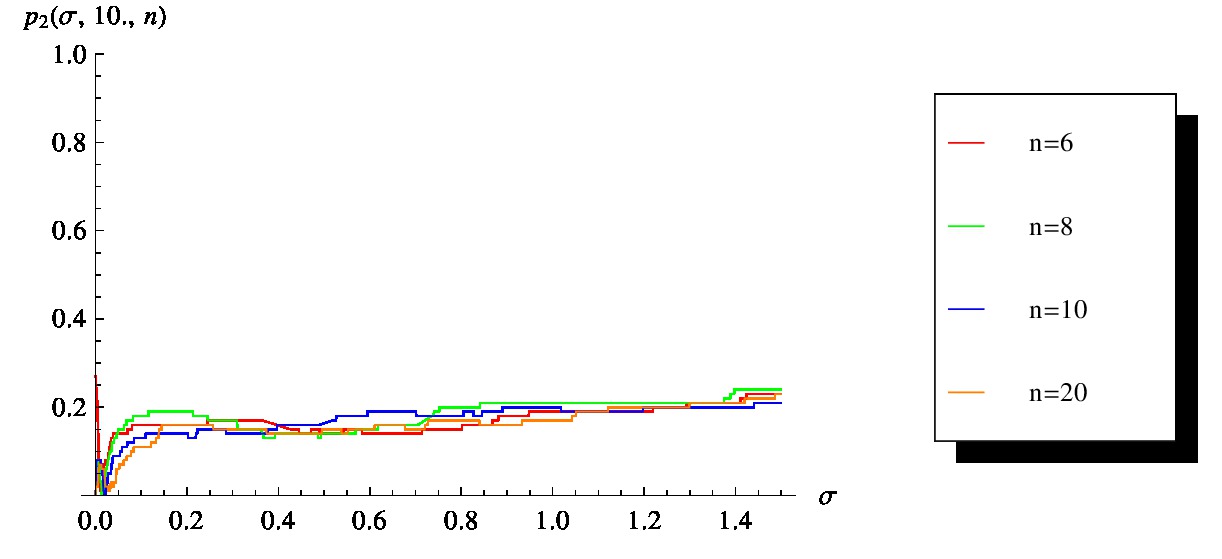}
\end{center}

\begin{center}
\includegraphics{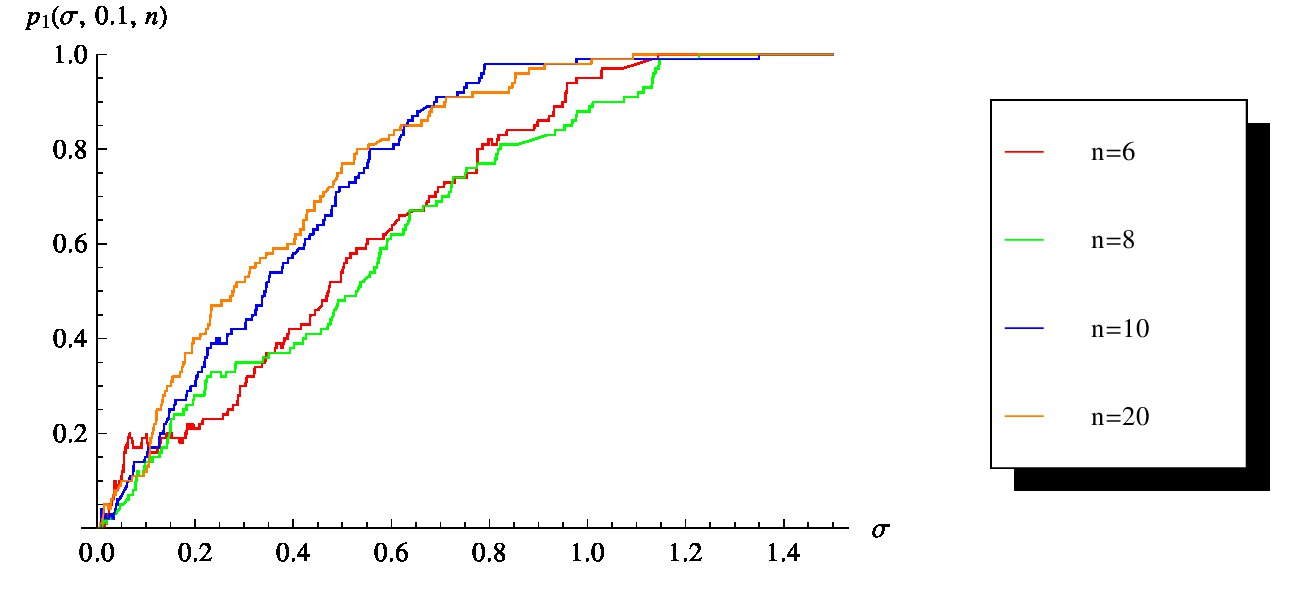}
\end{center}

\begin{center}
\includegraphics{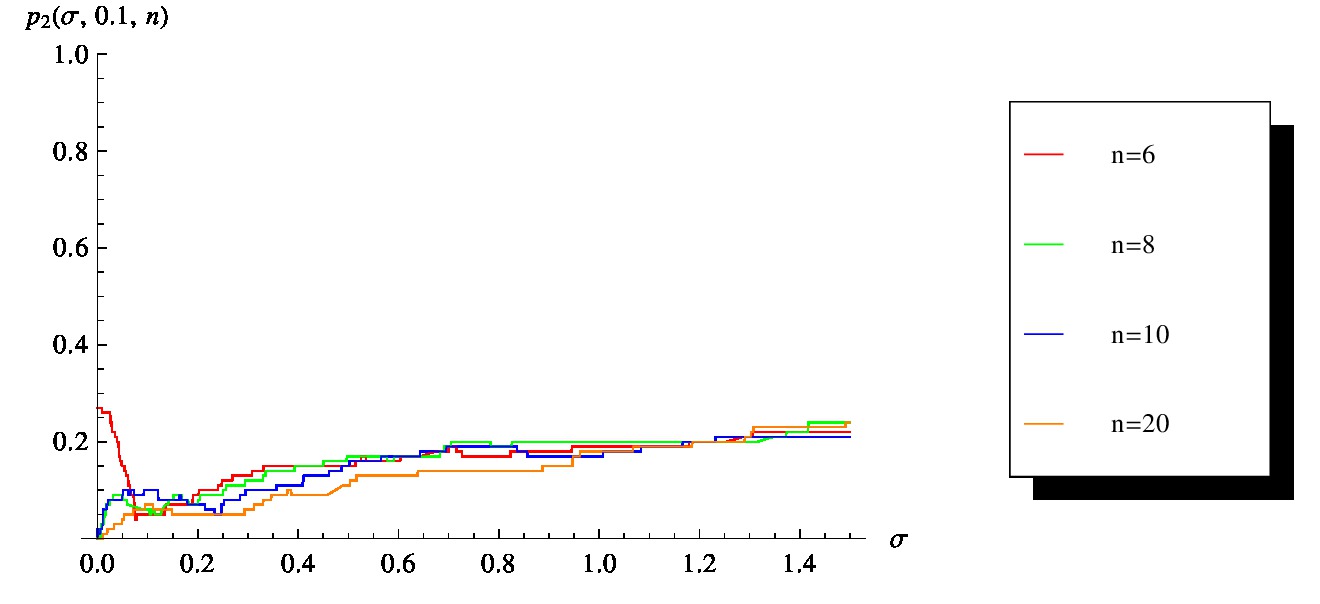}
\end{center}

\noindent  $\bullet$ Plot of $p_3(\sigma,n)$ and $p_4(\sigma,n).$

\begin{center}
\includegraphics{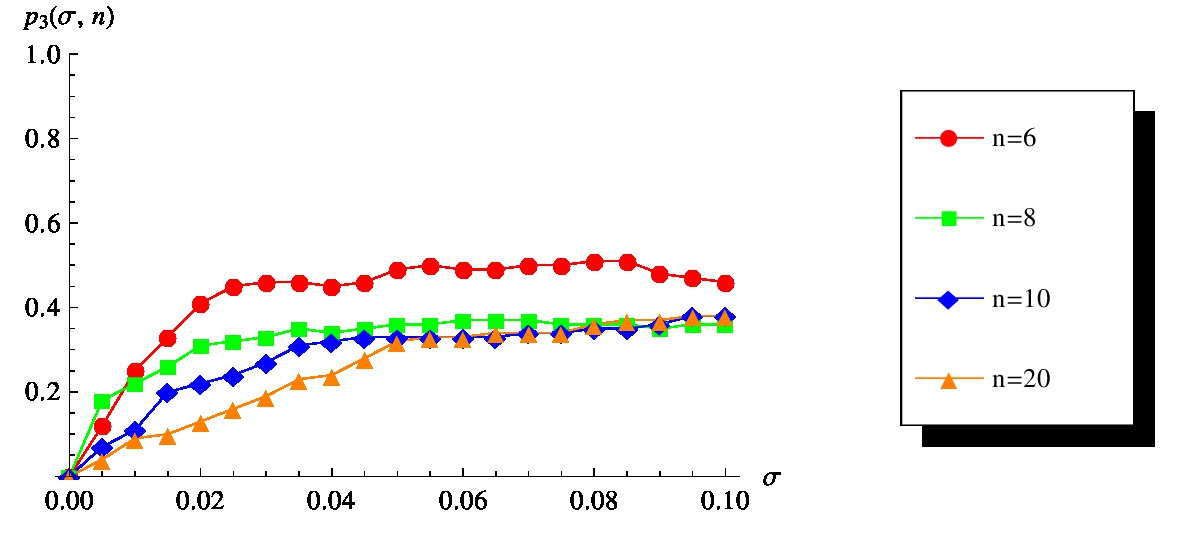}
\end{center}

\begin{center}
\includegraphics{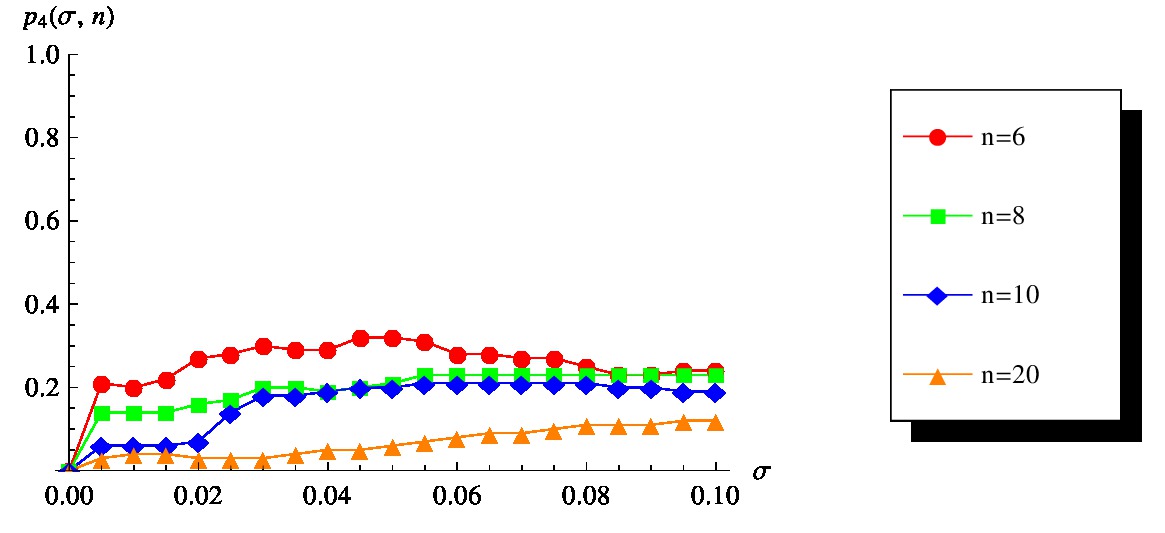}
\end{center}

\end{document}